\documentclass[12pt]{article}
\usepackage{latexsym,a4wide,amssymb,amsmath}

\newtheorem{thm}{Theorem}
\newtheorem{lem}{Lemma}
\newtheorem{cor}{Corollary}
\newtheorem{conj}{Conjecture}
\newtheorem{prop}{Proposition}
\newtheorem{exer}{Exercise}
\newcommand{\ebox}{\hfill $\Box$\\\vspace{0.15cm}}
\newcommand{\pr}{{\bf Proof.}\ }
\newcommand{\bt}{\begin{thm}}
\newcommand{\et}{\end{thm}}
\newcommand{\bl}{\begin{lem}}
\newcommand{\el}{\end{lem}}
\newcommand{\bp}{\begin{prop}}
\newcommand{\ep}{\end{prop}}
\newcommand{\bc}{\begin{cor}}
\newcommand{\ec}{\end{cor}}
\newcommand{\bcj}{\begin{conj}}
\newcommand{\ecj}{\end{conj}}
\newcommand{\bex}{\begin{exer}}
\newcommand{\eex}{\end{exer}}
\newcommand{\bi}{\begin{itemize}}
\newcommand{\ei}{\end{itemize}}
\newcommand{\be}{\begin{equation}}
\newcommand{\ee}{\end{equation}}
\newcommand{\ben}{\begin{enumerate}}
\newcommand{\een}{\end{enumerate}}

\newcommand{\mt}{t\kern-0.035cm\char39\kern-0.03cm}
\newcommand{\ml}{l\kern-0.035cm\char39\kern-0.03cm}
\newcommand{\md}{d\kern-0.035cm\char39\kern-0.03cm}

\newcommand{\noi}{\noindent}

\begin{document}

\title{\vspace{-2.3cm} Asymptotically approaching the Moore bound \\ for diameter three by Cayley graphs}

\author{}
\date{}
\maketitle

\begin{center}
\vspace{-1.3cm}

{\large Martin Bachrat\'y} \\
\vspace{2mm} {\small Comenius University, Bratislava, Slovakia}\\

\vspace{4mm}

{\large Jana \v Siagiov\'a} \\
\vspace{2mm} {\small Slovak University of Technology, Bratislava, Slovakia}

\vspace{4mm}

{\large Jozef \v Sir\'a\v n} \\
\vspace{2mm} {\small
Open University, Milton Keynes, U.K., and \\ Slovak University of Technology, Bratislava, Slovakia}

\vspace{4mm}

\end{center}

\begin{abstract}

The largest order $n(d,k)$ of a graph of maximum degree $d$ and diameter $k$ cannot exceed the Moore bound, which has the form $M(d,k)=d^k - O(d^{k-1})$ for $d\to\infty$ and any fixed $k$. Known results in finite geometries on generalised $(k+1)$-gons imply, for $k=2,3,5$, the existence of an infinite sequence of values of $d$ such that $n(d,k)=d^k - o(d^k)$. This shows that for $k=2,3,5$ the Moore bound can be asymptotically approached in the sense that $n(d,k)/M(d,k)\to 1$ as $d\to\infty$;  moreover, no such result is known for any other value of $k\ge 2$. The corresponding graphs are, however, far from vertex-transitive, and there appears to be no obvious way to extend them to vertex-transitive graphs giving the same type of asymptotic result.

The second and the third author (2012) proved by a direct construction that the Moore bound for diameter $k=2$ can be asymptotically approached by Cayley graphs. Subsequently, the first and the third author (2015) showed that the same construction can be derived from generalised triangles with polarity.

By a detailed analysis of regular orbits of suitable groups of automorphisms of graphs arising from polarity quotients of incidence graphs of generalised quadrangles with polarity, we prove that for an infinite set of values of $d$ there exist Cayley graphs of degree $d$, diameter $3$, and order $d^3{-}O(d^{2.5})$. The Moore bound for diameter $3$ can thus as well be asymptotically approached by Cayley graphs. We also show that this method does not extend to constructing Cayley graphs of diameter $5$ from generalised hexagons with polarity.

\vskip 3mm

\noi {\em Keywords:} Degree; Diameter; Moore bound; Cayley graph; Generalised quadrangle; Automorphism; Polarity.

\end{abstract}

\vskip 3mm

\section{Introduction}\label{intro}

For positive integers $d$ and $k$ let $n(d,k)$ denote the largest order of a graph of maximum degree $d$ and diameter $k$. It is well known that the value of $n(d,k)$ cannot exceed the Moore bound $M(d,k)=1+d+d(d-1)+\ldots +d(d-1)^{k-1}$. Setting trivial cases aside, for $d\ge 3$ and $k\ge 2$ we have $n(d,k)=M(d,k)$ only for $k=2$ and $d=3$, $7$, and possibly $57$, the unique graphs for the first two degrees being the Petersen and the Hoffman-Singleton graph \cite{HS,BI,Da}. For a survey of results about (im)possibility of getting `close' to the Moore bound for the remaining values of $d$ and $k$ we refer to \cite{MS}. The main driving forces in this field appear to be the question of Bermond and Bollob\'as \cite{BB} if for arbitrarily large $c$ there exist $d,k$ such that $n(d,k)<M(d,k)-c$, and the problem of Delorme \cite{De} of determining the value of ${\rm lim\ sup}_{d\to\infty}\ n(d,k)/d^k={\rm lim\ sup}_{d\to\infty}\ n(d,k)/M(d,k)$ for every fixed $k\ge 2$.
\smallskip

Regarding the question of Bermond and Bollob\'as, a substantial progress has recently been made in \cite{E+} by proving that for any fixed $d$ and any $c>0$ the order of the largest {\sl vertex-transitive} $d$-regular graph of diameter $k$ is smaller than $M(d,k)-c$ for almost all $k$. The best available result addressing Delorme's problem is his own observation \cite{De} that ${\rm lim\ sup}_{d\to\infty}\ n(d,k)/M(d,k) = 1$ for $k\in \{2,3,5\}$. This follows by taking polarity quotients of the incidence graphs of generalised $n$-gons admitting a polarity (cf. \cite{M}) for $n\in\{3,4,6\}$, respectively. Despite having a fairly large automorphism group compared to their order, these graphs are not even regular and by \cite{BS} there appears to be no obvious way to extend them to vertex-transitive graphs by just adding edges in the case of diameter $k=2$.
\smallskip

For the remaining diameters the best currently known results on Delorme's problem are much weaker but far from easy to prove. We know that ${\rm lim\ sup}_{d\to\infty}\ n(d,4)/M(d,4) \ge 1/4$ by \cite{De2},  and from \cite{CG} we have ${\rm lim\ sup}_{d\to\infty}\ n(d,k)/M(d,k) \ge (1.6)^{-k}$ for $k\ge 6$, where $1.6$ can be replaced by $1.57$ for $k\equiv -1,0,1$ mod $6$.
\smallskip

In the light of the above-mentioned result of \cite{E+} addressing the question of Bermond and Bollob\'as for vertex-transitive graphs, it is natural to ask if the Moore bound can be {\sl asymptotically approached}, in the sense of Delorme's limit superior being equal to $1$, by {\sl vertex-transitive}, or even {\sl Cayley} graphs. The importance of this direction of research is underscored by the fact that, from the practical point of view, computer generation of record large graphs of given degree and diameter is almost exclusively limited to searching over Cayley graphs in cases when the degree or diameter are beyond values manageable by other methods; cf. \cite{LoSi} and the on-line tables \cite{table}.
\smallskip

To this end, for $d\ge 3$ and $k\ge 2$ we let $vt(d,k)$ and $Cay(d,k)$ denote the largest order of a vertex-transitive and a Cayley graph, respectively, of degree $d$ and diameter $k$; clearly, $vt(d,k)\ge Cay(d,k)$. The task now is to estimate the values of ${\rm lim\ sup}_{d\to\infty}\ vt(d,k)/M(d,k)$ and ${\rm lim\ sup}_{d\to\infty}\ Cay(d,k)/M(d,k)$ for $k\ge 2$ and, specifically, to determine if the Moore bound can be asymptotically approached by vertex-transitive or Cayley graphs for diameters $2$, $3$ and $5$.  Here, however, the available results are scarcer and, expectedly, not as good as those for $n(d,k)$. Let us begin with $k\ge 3$. In the vertex-transitive case, the digraphs of \cite{F+} yield, after ignoring edge directions, ${\rm lim}_{d\to\infty}\ vt(d,k)/M(d,k)\ge 2^{-k}$ for every $k\ge 3$. For Cayley graphs, constructions of \cite{M+,M++} give ${\rm lim}_{d\to\infty}\ Cay(d,k)/M(d,k)\ge k\cdot 3^{-k}$ for every $k\ge 3$, with improvements of the lower bounds by \cite{V} to $3\cdot 2^{-4}$, $32\cdot 5^{-4}$ and $25\cdot 4^{-5}$ for $k=3$, $4$ and $5$, respectively.
\smallskip

For $k=2$ the strongest finding in this category is ${\rm lim\ sup}_{d\to\infty}\ Cay(d,2)/M(d,2)=1$, showing that the Moore bound for diameter $2$ can be asymptotically approached by Cayley graphs. This result was obtained in \cite{Sia} by a direct construction of certain Cayley graphs of one-dimensional affine groups over finite fields of characteristic $2$. Later in \cite{BS} it was shown that the construction of \cite{Sia} is equivalent to extending a regular orbit of a polarity quotient of the incidence graph of a generalised triangle under the action of a suitable group.
\smallskip

Our aim is to show that the Moore bound for diameter $3$ can also be asymptotically approached by Cayley graphs, i.e., to prove that ${\rm lim\ sup}_{d\to\infty}\ Cay(d,3)/M(d,3) = 1$. In fact, we prove that for an infinite set of values of $d$ there exist Cayley graphs of degree $d$, diameter $3$, and order $d^3-O(d^{2.5})$. The method is a variant of the one used in \cite{BS}, namely, extension of a regular orbit of a suitable subgroup of the automorphism group of a polarity quotient of the incidence graph of a generalised quadrangle. Details, however, are much more subtle and complex in comparison with those of \cite{BS}. We also show that an extension of this method is not feasible for proving an analogous result for diameter $5$.
\smallskip

The paper is organised as follows. In Section 2 we review basic concepts on the finite generalised quadrangles with polarity and their automorphisms. These are used in Section 3 to study incidence in  auxiliary graphs obtained from the incidence graphs of finite generalised quadrangles by polarity. In Section 4 we investigate induced subgraphs of the auxiliary graphs obtained as orbits of suitable groups of automorphisms. The induced subgraphs are finally extended in Section 5 to give Cayley graphs which prove that ${\rm lim\ sup}_{d\to\infty}\ Cay(d,3)/M(d,3) = 1$. We conclude by showing that applying this scenario to generalised hexagons with polarity does not produce Cayley graphs that would asymptotically approach the Moore bound for diameter $5$.
\smallskip

\section{Generalised quadrangles, polarity and symmetries}\label{prelim}

Let $q$ be a prime power and let $F=GF(q)$ be the Galois field of order $q$. As usual, let $F^+$ and $F^*$ be the additive and the multiplicative group of $F$. We begin by recalling the projective geometry $PG(3,q)$ whose points are the $1$-dimensional subspaces of $F^4$ minus the origin (sometimes called projective vectors), that is, equivalence classes ${\bf [x]}$ of non-zero quadruples ${\bf x}=(x_0,x_1,x_3,x_3)\in F^4$, with two quadruples ${\bf x}$ and ${\bf y}$ equivalent if $y_i=tx_i$ for some $t\in F^*$ and every $i\in \{0,1,2,3\}$. A subset $S$ of $PG(3,q)$ is {\em totally isotropic} if for any two points ${\bf [x]}, {\bf [y]} \in S$ we have $x_0y_1-x_1y_0+x_2y_3-x_3y_2=0$. Note that each of the $q^3+q^2+q+1$ points are totally isotropic themselves, and an easy counting argument shows that there are exactly $q^3+q^2+q+1$ totally isotropic lines (that is, $2$-dimensional subspaces of $F^4$ with the origin removed) of $PG(3,q)$.
\smallskip

Total isotropy helps us introduce an important incidence geometry within $PG(3,q)$, standardly denoted $W(q)$; see e.g. \cite{M,PT}.  Points of $W(q)$ are the points of $PG(3,q)$, lines of $W(q)$ are the totally isotropic lines of $PG(3,q)$, and incidence is defined by containment as in $PG(3,q)$. Every line of $W(q)$ contains $q+1$ points (as in $PG(3,q)$) and, by counting, every point of $W(q)$ lies on $q+1$ lines of $W(q)$. The incidence structures $W(q)$ are prominent examples of {\em generalised quadrangles} \cite{M,PT}.
\smallskip

A polarity $\pi$ of $W(q)$ is an involutory mapping that sends the point set of $W(q)$ onto its line set and vice versa, with the property that for any two points $u$, $v$ lying on a line $\ell$ of $W(q)$, the lines $\pi(u)$ and $\pi(v)$ intersect at the point $\pi(\ell)$. By a classical result of Tits \cite{Ti}, the incidence structure $W(q)$ admits a polarity if and only if $q$ is an odd power of $2$. We describe such a polarity next, following \cite{Pa}.
\smallskip

From now on and throughout the paper, let $q=2^{2n+1}$ for some positive integer $n$. Let $\omega=2^{n+1}$ and let $\sigma$ be the automorphism of $F=GF(q)$ given by $\sigma(x)=x^{\omega}$, so that $\sigma^2(x)=x^2$. For every point $u={\bf [x]}\in W(q)$ let $c=x_0x_1+x_2x_3$ and let $\pi(u)$ be the set of all non-zero vectors of $F^4$ spanned by the totally isotropic set of four vectors \begin{equation}\label{span} (0,c^{\omega/2},x_0^{\omega},x_2^{\omega}),\ (c^{\omega/2},0,x_3^{\omega},x_1^{\omega}),\  (x_0^{\omega},x_3^{\omega},0,c^{\omega/2})\ {\rm and}\ (x_2^{\omega},x_1^{\omega},c^{\omega/2},0)\ ;\end{equation} one may check that this is indeed a line of $W(q)$. In the reverse direction, let $\ell$ be a line of $W(q)$ through a pair of distinct points ${\bf [x]}$ and ${\bf [y]}$, and let $\delta_{ij}=x_iy_j+x_jy_i$ for any distinct $i,j\in \{0,1,2,3\}$. Then,  $\pi(\ell)$ is the point ${\bf [z]}\in W(q)$ such that \begin{equation}\label{pol} z_0=\delta_{02}^{\omega/2},\  z_1=\delta_{31}^{\omega/2}, \ z_2=\delta_{03}^{\omega/2}\ {\rm and}\  z_3=\delta_{21}^{\omega/2}\ .\end{equation} It can be shown (see \cite{Pa} for a large number of details) that the mapping $\pi$ is a polarity of $W(q)$.
\smallskip

For every $x,y\in F$ we let $f(x,y)=x^{\omega+2}+xy+y^{\omega}$. The set of matrices $M(r;a,b)$ given, for all $r\in F^*$ and $a,b\in F$, by

\begin{equation}\label{M}
M(r;a,b)=\left( \begin{array} {cccc}
1 & f(a,b) & a & b \\
0 & r^{\omega+2} & 0 & 0 \\
0 & (a^{\omega+1}{+}b)r & r & a^{\omega}r \\
0 & ar^{\omega+1} & 0 & r^{\omega+1}
\end{array}\right)
\end{equation}

\noi is closed under multiplication and forms a group $G$ of order $q^2(q-1)$. In particular, one can verify that \[ M(r;a,b)M(s;c,d)=M(rs;as+c,bs^{\omega+1}+d+ac^{\omega}s) \] and so $G$ is isomorphic to an iterated split extension of the form $(F^+ \rtimes F^+)\rtimes F^*$. The group $G$ acts on $W(q)$ as a group of collineation by right multiplication. In \cite{Ti} J. Tits proved that the group of all collineations of $GP(3,q)$ leaving the set $\Omega=\{[0,1,0,0]\}\cup\{[1,f(x,y),x,y];\ x,y\in F\}$  invariant is (isomorphic to) the Suzuki group $Sz(q)={^2}B_2(q)$, a simple group of order $q^2(q^2+1)(q-1)$. Moreover, in the above representation, $G$ is the subgroup of $Sz(q)$ stabilising the point $[0,1,0,0]$. We note that $\Omega$ is the set of absolute points with respect to $\pi$, that is,  the set of points $u$ for which $u\in \pi(u)$; it is also known as the Suzuki-Tits ovoid.
\smallskip

The action of $G$ on the points (and hence also on the lines) of $W(q)$ is intransitive and a straightforward calculation shows that $G$ has the following five orbits $O_1$ - $O_5$ on points of $W(q)$:
\begin{itemize}
\item $O_1=\{[1,f(x,y),x,y];\ u,v\in F\}=\Omega{\setminus}\{[0,1,0,0]\}$, of size $q^2$;
\item $O_2=\{[0,x,1,y];\ x,y\in F\}$, of size $q^2$;
\item $O_3=\{[0,x,0,1];\ x\in F\}$, of size $q$;
\item $O_4=\{[0,1,0,0]\}$, which is the unique fixed point of $G$; and
\item $O_5=V \setminus (\cup_{i=1}^{4}O_i)$, of size $q^2(q-1)$.
\end{itemize}
We point out the (for us) important fact that $G$ acts regularly on the orbit $O_5$.
\smallskip

\section{The graph arising from factorisation by polarity}

Let $A(q)$ be the graph whose vertex set $V$ is the set of all points $u$ of $W(q)$, with two distinct vertices $u,v$ adjacent in $A(q)$ if $u\in \pi(v)$, which, by properties of $\pi$ is equivalent to $v\in \pi(u)$. Observe that every vertex of $A(q)$ is adjacent to $q+1$ or $q$ vertices, and the degree of a vertex $u$ in $A(q)$ is $q$ if and only if $u\in \pi(u)$, that is, if and only if $u$ is an element of the Suzuki-Tits ovoid, of cardinality $q^2+1$. Note that we could have defined $A(q)$ as the quotient graph of the bipartite point-line incidence graph of $W(q)$ obtained by factorisation by the polarity $\pi$, that is, by identifying $u$ with $\pi(u)$ throughout and suppressing eventual edges between $u$ and $\pi(u)$.
\smallskip

A basic property of the generalised quadrangle $W(q)$ (cf. \cite{PT}) is that for any line $\ell$ and any point $u$ not on $\ell$ there is a unique point $u'\in \ell$ such that $u$ and $u'$ are collinear in $W(q)$. This immediately translates to the observation that the diameter of the graph $A(q)$ is equal to $3$. Indeed, suppose that $u,v\in V$ are vertices that are not adjacent in the graph $A(q)$. Then, since the line $\ell=\pi(v)$ does not contain $u$, by the above property there is a (unique, which is not important for this argument) point $u'\in \ell$ such that $u,u'\in \ell'$ for some line $\ell'\in W(q)$. But letting $\pi(\ell')=v'$ we have a path $uv'u'v$ of length $3$ in $A(q)$. (Note that, in the above argument, the vertices $u$ and $v$ might correspond to collinear points of $W(q)$, in which case we would have $u$ and $v$ joined by a path of length $2$ and also by a path of length $3$, giving rise to a cycle of length $5$ in $A(q)$.)
\smallskip

For further analysis we will need a description of the neighbourhood $N(u)$ of a few vertices $u$ of $A(q)$, that is, the set of all $v\in V$ adjacent to $u$ in $A(q)$. By the adjacency rule in $A(q)$, a vertex $u=[{\bf x}]$ is adjacent to precisely the vertices $v=[{\bf y}]\ne [{\bf x}]$ that are points on the line $\pi[{\bf x}]$ spanned by the vectors in (\ref{span}). Informally, a projective vector $[{\bf x}]$ is adjacent in $A(q)$ precisely to the projective vectors in $\pi[{\bf x}]$ distinct from ${\bf x}$. A concrete identification of $\pi[{\bf x}]$ can be cumbersome in general but is not for the few vertices we need. We illustrate the process on the vertex $[{\bf x}]=[1,1,1,1]$. Taking the first and the third vector in (\ref{span}) one sees that $\pi[{\bf x}]$ is generated by the two (totally isotropic) projective vectors $[0,0,1,1]$ and $[1,1,0,0]$. Excluding the self-adjacency it follows that  the vertex $[1,1,1,1]$ is adjacent in $A(q)$ to the $q$ vertices $[0,0,1,1]$ and $[1,1,z,z]$ for $z\in F{\setminus}\{1\}$, or, equivalently, $N[1,1,1,1]=\{[0,0,1,1]\}\cup\{[1,1,z{+}1,z{+}1];\ z\in F^*\}$.
In a completely analogous way we obtain the following table.

\begin{center}
\begin{tabular}{|c|c|}
\hline
$ u\in V$ & $N(u)$ \\
\hline
$[1,0,0,0]$ & $\{[0,0,1,0]\}\cup\{[1,0,z,0];\ z\in F^*\}$ \\
\hline
$[0,1,0,0]$ & $\{[0,z,0,1];\ z\in F\}$ \\
\hline
$[0,0,1,0]$ & $\{[0,0,0,1]\}\cup\{[1,0,0,z];\ z\in F\}$ \\
\hline
$[0,0,0,1]$ & $\{[0,1,0,0]\}\cup\{[0,z,1,0];\ z\in F\}$ \\
\hline
$[1,1,0,0]$ & $\{[0,1,1,0]\}\cup\{[1,z,z,1];\ z\in F\}$ \\
\hline
$[0,1,1,0]$ & $\{[0,0,0,1]\}\cup\{[1,1,0,z];\ z\in F\}$ \\
\hline
$[0,0,1,1]$ & $\{[0,1,0,1]\}\cup\{[1,z,1,z];\ z\in F\}$ \\
\hline
$[1,1,0,1]$ & $\{[0,1,1,0]\}\cup\{[1,1{+}z,z,1];\ z\in F^*\}$ \\
\hline
$[1,1,1,1]$ & $\{[0,0,1,1]\}\cup\{[1,1,z{+}1,z{+}1];\ z\in F^*\}$ \\
\hline
\end{tabular}
\end{center}

\centerline{Table 1: Neighbourhood of selected vertices of $A(q)$.} \medskip

As the next step we show that the group $G$ generated by the matrices $M=M(r;a,b)$ introduced in (\ref{M}) acts on the vertices of $A(q)$ by right multiplication. This is clearly equivalent to the statement that $\pi[{\bf x}M]=\pi[{\bf x}]M$ for every $[{\bf x}]\in V$. To verify this it is sufficient to restrict ourselves to a particular $[{\bf z}]$ from each orbit $O\in \{O_1,\ldots,O_5\}$ listed in section \ref{prelim}. Indeed, if we choose a $[{\bf z}]\in O$, then, for every $[{\bf x}]\in O$ we have $[{\bf x}]=[{\bf z}M_z]$ for some $M_z\in G$. Provided we show that $\pi[{\bf z}M']=\pi[{\bf z}]M'$ for every $M'\in G$, we then have $\pi[{\bf x}M]= \pi[{\bf z}M_zM] = \pi[{\bf z}]M_zM = \pi[{\bf z}M_z]M = \pi[{\bf x}]M$ for every $M\in G$. The verification can now be done by letting $[{\bf z}]$ be the representatives $[1,0,0,0]$, $[0,0,1,0]$, $[0,0,0,1]$, $[0,1,0,0]$ and $[1,1,0,0]$ of the orbits $O_1$ to $O_5$, respectively.
\smallskip

We illustrate the procedure on the computationally most demanding case when $[{\bf z}]=[1,1,0,0]$. Instead of proving $\pi[{\bf z}M]=\pi[{\bf z}]M$ we prove the equivalent equality $[{\bf z}]M=\pi[\pi({\bf z})M]$. Looking at Table 1 we see that both vertices $[0,1,1,0]$ and $[1,0,0,1]$ are neighbours of $[{\bf z}]=[1,1,0,0]$ in $A(q)$. This is equivalent to stating that $\pi[{\bf z}]$ is the line of $W(q)$ through the points $[{\bf x}]=[0,1,1,0]$ and $[{\bf y}]=[1,0,0,1]$. Since $G$ is a collineation group of $W(q)$, it follows that for every $M=M(r;a,b)$ from (\ref{M}) the point $[{\bf z}]M=[1,f(a,b)+r^{\omega+2},a,b]$ lies on the line $\ell=\pi[{\bf z}]M$ through the points $[{\bf x}]M=[0,r^{\omega+2} + r(a^{\omega+1}+b), r, a^{\omega}r]$ and $[{\bf y}]M=[1,f(a,b)+ ar^{\omega+1}, a,b+r^{\omega+1}]$. We now calculate the point $\pi(\ell)$ by the procedure described in (\ref{pol}). A somewhat lengthy but straightforward verification shows that $\delta_{02}^{\omega/2}=r^{\omega/2}$, $\delta_{31}^{\omega/2} = r^{\omega/2}(f(a,b)+r^{\omega+2})$, $\delta_{03}^{\omega/2}=r^{\omega/2}a$ and $\delta_{21}^{\omega/2}= r^{\omega/2}b$, so that $[ \delta_{02}^{\omega/2}, \delta_{31}^{\omega/2}, \delta_{03}^{\omega/2}, \delta_{21}^{\omega/2}] = [1,f(a,b)+r^{\omega+2},a,b]$. Thus, $\pi(\ell)=[{\bf z}]M$ and hence $\pi[\pi({\bf z})M]=[{\bf z}]M$, which is what we wanted to establish. We leave out the details for the remaining choices of orbit representatives as they are similar (and easier).
\smallskip

\section{The subgraph induced by the regular orbit}

From this point on we will focus on the subgraph $B(q)$ of $A(q)$ induced by the subset $O_5$ of $V$. Since the group $G$ acts regularly on the vertex set of $B(q)$, it follows \cite{Sab} that $B(q)$ is a Cayley graph for the group $G$ and some generating set $S$ for $G$. The set $S$ is uniquely determined and can be recovered by looking at the `local' action of $G$ on a vertex as follows. 
\smallskip

Let $u=[{\bf x}]$ be a fixed vertex of $B(q)$; for later convenience we will fix the vertex $u=[1,1,0,0]$ of $B(q)$ throughout. By the regular action of $G$ on vertices of $B(q)$, for every neighbour $w$ of $u$ in $B(q)$ there is exactly one $g_w\in G$ represented by a matrix $M_w$ as in (\ref{M}) such that $w=ug_w=[{\bf x}]M_w$. To see what happens with adjacency in an arbitrary vertex, just apply an arbitrary element $g\in G$ to this situation to conclude that a vertex $ug$ is adjacent to $wg=u(g_wg)$ for every neighbour $w$ of $u$. Let now $S$ be the set of all the $g_w\in G$ where $w$ ranges over all neighbours of $u$ in $B(q)$. The regular action of $G$ enables us to identify $G$ with the vertex set of $B(q)$ by means of the bijection $g\mapsto ug$ and the same bijection gives an isomorphism of the Cayley graph $C(G,S)$ onto $B(q)$; note that this isomorphism maps the identity of $G$ onto the fixed vertex $u$. The adjacency rule in $C(G,S)$ follows from the above, namely, $g\in G$ is adjacent to $g_wg$ for every $g_w\in S$.
\smallskip

By Table 1 and the description of the orbit $O_5$, the neighbourhood of $u$ in $B(q)$ is the set $\{[1,x,x,1];\ x\in F,\ x\ne 1\}$. The elements $g_w\in S$ are thus matrices $M(r;a,b)$ from (\ref{M}) such that $[1,1,0,0]M(r;a,b)=[1,x,x,1]$, $x\ne 1$. A quick calculation reveals that we must have $b=1$, and $a,r$ are tied by the equation $a^{\omega+2}+r^{\omega+2}=1$. This equation has a unique solution $a=a(r)\ne 1$ for every $r\in F^*$, given by $a(r)=(1+r^{\omega+2})^{1-\omega/2}$ (recall that every element in $F^*$ has a unique square root). One can verify that $[1,1,0,0]M(r;a(r),1)=[1,a,a,1]$, $a\ne 1$. We have therefore identified the generating set as $S=\{M(r;a(r),1);\ r\in F^*\}$; this set is closed under taking inverses by the above construction. Summing up, we have:

\bl\label{cay}
The graph $B(q)$ is isomorphic to the Cayley graphs $C(G,S)$ of degree $q-1$ with the generating set $S=\{M(r;a(r),1);\ r\in F^*\}$. \hfill $\Box$
\el

In contrast to $A(q)$ the diameter of $B(q)$ is larger than $3$. This is a consequence of the following fact. Assume that $v$ and $w$ are vertices of $A(q)$ such that $v$ represents a point of $W(q)$ and $w$ a line of $W(q)$ not containing the point and such that $v$ and $\pi(w)$ are not on a line of $W(q)$. Then, by properties of generalised quadrangles, there is a unique path from $v$ to $w$ of length $3$ in $A(q)$, and no shorter path, and removing vertices from $A(q)$ will destroy some of these paths. Our aim is to extend the generating set $S$ of the Cayley graph $C(G,S)\simeq B(q)$ by just a `few' new generators to ensure that the new Cayley graph (which will contain $B(q)$ as a spanning subgraph) has diameter $3$. We thus need to identify pairs of vertices of $B(q)$ that could end up at distance larger than $3$ after the removal of the sets $O_1$ -- $O_4$ from $A(q)$. Since $B(q)$ is a Cayley graph, it is sufficient to do this for pairs of vertices $u,v\in O_5$ in which $u=[1,1,0,0]$ is our fixed vertex.
\smallskip

In the considerations below we will frequently refer to Table 1 without explicit alerts. In the graph $A(q)$, our fixed vertex $u=[1,1,0,0]\in O_5$ is adjacent to the two vertices $u_1=[1,1,1,1]\in O_1$, $u_2=[0,1,1,0]\in O_2$, and has another $q-1$ neighbours inside $B(q)$, all of the form $[1,x,x,1]$ for $x\in F$, $x\ne 1$. In particular, the degree of the Cayley graph $B(q)$ is $q-1$. The neighbour $u_1=[1,1,1,1]\in O_1$ of $u$ in $A(q)$ is adjacent to $u_2'=[0,0,1,1]\in O_2$ and to another $q-1$ vertices $[1,1,z{+}1,z{+}1]$, $z\in F^*$, in $B(q)$. It follows that when restricting to $B(q)$ we lose the paths of length $2$ from $u$ via $u_1$ to the $q-2$ vertices in the set $L(u_1)=N(u_1){\setminus}\{u,u_2'\} = \{[1,1,z{+}1,z{+}1];\ z\in F^*{\setminus}\{1\}\}$. Similarly, the neighbour $u_2\in O_2$ of $u$ in $A(q)$ is adjacent to $u_3= [0,0,0,1]\in O_3$, $u'_1=[1101]\in O_1$, and has further $q-1$ neighbours $[1,1,0,z{+}1]$, $z\in F^*$, in $B(q)$. Again, when restricting to $B(q)$ we lose the paths of length $2$ from $u$ through $u_2$ to the $q-2$ vertices in the set $L(u_2)=N(u_2){\setminus} \{u,u_3,u'_1\}= \{[1,1,0,z{+}1];\ z\in F^*{\setminus}\{1\}\}$.
\smallskip

The vertex $u'_1$ is adjacent in $A(q)$ to $u_2=[0,1,1,0]\in O_2$ and to the $q-1$ vertices forming the set $L(u'_1)= N(u'_1) {\setminus} \{u_2\} = \{[1,z{+}1,z,1];\ z\in F^*\}$. The neighbourhood of $u'_2$ in $A(q)$ consists of $u_1=[1,1,1,1]\in O_1$, $u'_3=[0,1,0,1]\in O_3$, and the $q-1$ vertices in the set $L(u'_2)=N(u'_2){\setminus} \{u_1,u'_3\} =\{[1,z{+}1,1,z{+}1];\ z\in F^*\}$. This implies that when restricting our attention to the graph $B(q)$ we lose the paths of length $3$ of the form $uu_1u'_2w$ for $w\in L(u'_2)$ and $uu_2u'_1w$ for $w\in L(u'_1)$.
\smallskip

We have identified four sets of vertices, namely, $L(u_1)$, $L(u_2)$, $L(u'_1)$ and $L(u'_2)$, to which we lose paths from $u$ through $u_1$ and $u_2$ when considering $B(q)$. In order to make up for this we will suitably extend the generating set $S$ in the Cayley graph $C(G,S)\simeq B(q)$. To do so we first determine the action of $G$ on the four sets, which is equivalent to determining the vertex stabilisers in $G$ of $u_1$, $u_2$, $u'_1$ and $u'_2$.

\bl\label{stab}
The stabilisers in $G$ of $u_1$, $u_2$, $u'_1$ and $u'_2$ are all cyclic and isomorphic to $F^*$. In more detail:
\begin{itemize}
\item $Stab_G(u_1)=Stab_G(u'_2)=\{M(r;r{+}1,r{+}1);\ r\in F^*\}\simeq F^*$, acting regularly on both $L(u_1)\cup\{u\}$ and $L(u'_2)$  by $[1,1,z{+}1,z{+}1]\cdot M(r;r{+}1,r{+}1) = [1,1,zr{+}1,zr{+}1]$ and  $[1,z{+}1,1,z{+}1]\cdot M(r;r{+}1,r{+}1) = [1,zr^{\omega+1}+1,1,zr^{\omega+1}+1]$.
\item $Stab_G(u_2)=Stab_G(u'_1)=\{M(r;0,1{+}r^{\omega+1});\ r\in F^*\}\simeq F^*$, acting regularly on \linebreak both sets $L(u_2)\cup\{u\}$ and $L(u'_1)$ by $[1,1,0,z{+}1]\cdot M(r;0,1{+}r^{\omega+1}) = [1,1,0,zr^{\omega+1}]$ and \linebreak $[1,z{+}1,z,1]\cdot M(r;0,1{+}r^{\omega+1}) = [1,zr+1,zr,1]$ for $z\in F^*$.
\end{itemize}
\el

\pr Calculation of the stabilisers from the equation $[{\bf x}]M(r;a,b)=[{\bf x}]$ for $[{\bf x}]\in \{u_3,u_4,u'_3,u'_4\}$ and solving for $a,b$ in terms of $r$ is straightforward. Regarding isomorphism of the stabilisers with $F^*$, for $Stab_G(u_1)=Stab_G(u'_2)$ it is given by $\theta_1:\ r\mapsto M(r;r{+}1,r{+}1)$ since one can check that $M(r;r{+}1,r{+}1)M(s;s{+}1,s{+}1) =M(rs;rs{+}1,rs{+}1)$ for all $r,s\in F^*$. The same type of isomorphism $\theta_2:\ r\mapsto M(r;0,1{+}r^{\omega+1})$ works for $Stab_G(u_2)=Stab_G(u'_1)$, as one can verify that $M(r;0,1{+}r^{\omega+1})M(s;0,1{+}s^{\omega+1}) = M(rs;0,1{+}(rs)^{\omega+1})$ for every $r,s\in F^*$. The actions of the stabilisers on the four sets are obvious.
\ebox

\section{Large Cayley graphs of diameter 3 and small degree}

For any integer $m\ge 3$ let $c_m$ denote the smallest number of elements in a generating set $X$ of a Cayley graph of a cyclic group $C(Z_m,X)$ of diameter $2$. The best available general upper bound on $c_m$, which is, of course, the degree of $C(Z_m,X)$, is $c_m\le 2\lceil\sqrt{m}\rceil$, see \cite{Ers} for a short proof. Applying this to the cyclic group  $F^*$ of order $m=q-1$ we have the existence of a Cayley graph $C(F^*,X)$ of diameter $2$ with $|X|= c_{q-1} \le 2\lceil\sqrt{q-1}\rceil$. Taking images of $C(F^*,X)$ under the isomorphisms $\theta_1$ and $\theta_2$ from the proof of Lemma \ref{stab}, that is, letting $S_1=\theta_1(X)$ and $S_2=\theta_2(X)$ and denoting $H_1=Stab_G(u_1)= Stab_G(u'_2)$ and $H_2=Stab_G(u_2) =Stab_G(u'_1)$, we have constructed Cayley graphs $C(H_1,S_1)$ and $C(H_2,S_2)$ of degree $c_{q-1}$ and diameter $2$.
\smallskip

We are now ready to describe a suitable `small' extension of our generating set $S$ from Lemma \ref{cay} in the Cayley graph $C(G,S)\simeq B(q)$ to obtain a graph of diameter $3$. Let $M_1=M(1;1,1)$, $M_2=M(1;1,0)$, $S_3=\{M_1,M_1^{-1}, M_2,M_2^{-1}\}$ and let $S^*=S\cup S_1\cup S_2\cup S_3$; note that $S^*$ is closed under taking inverses.

\bt\label{main1}
For every $n\ge 1$ and $q=2^{2n+1}$ the Cayley graph $C(G,S^*)$ of order $q^2(q-1)$ and degree $q+2c_{q-1}+3$ has diameter $3$.
\et

\pr
The Cayley graph $C(G,S^*)$ has order $|G|=q^2(q-1)$, degree $|S|$ that evaluates to $q+2c_{q-1}+3$, and it contains the Cayley graph $C(G,S)\simeq B(q)$ as a spanning subgraph. We may thus identify the vertex set of $C(G,S^*)$ with the orbit $O_5$ of the group $G$ on the set $V$ as we did in the isomorphism from $C(G,S)$ onto $B(q)$.
\smallskip

We proceed by extending some of our earlier observations made for the fixed vertex $u$ to arbitrary vertices of $B(q)$. By the action of $G$ on vertices of the entire graph $A(q)$ and by the regular action of the same group on the vertex set $O_5$ of the graph $B(q)$, every vertex $v\in O_5$ has one neighbour $v_1\in O_1$, one neighbour $v_2\in O_2$, and $q-1$ neighbours within $B(q)$. Further, in $A(q)$, the vertex $v_1$ is adjacent to one vertex $v'_2\in O_2$ and to a set $L(v_1)$ of $q-2$ vertices in $B(q)$ distinct from $v$, and the vertex $v_2$ has one neighbour $v_3\in O_3$, another neighbour $v'_1\in O_1$ and $q-2$ neighbours in a set $L(v_2)$ of vertices of $B(q)$ distinct from $v$. Moreover, by our isomorphism from the Cayley graph $C(G,S)$ onto $B(q)$ constructed earlier, the subsets $L(v_1)$ and $L(v_2)$ of $O_5$ are $g$-images of the sets $L(u_1)$ and $L(u_2)$ for the element $g\in G$ that carries $u$ onto $v=ug$. Also, note that the unique vertex $[0,1,0,0]$ of $O_4$ of degree $q$ is incident only to $q$ vertices in $O_3$ and the neighbours of every vertex of $O_3$ (except $[0,1,0,0]$) lie in $O_2$.
\smallskip

Equipped with this we now prove that the diameter of $C(G,S^*)$ is $3$. We know that it is sufficient to check distances from our fixed vertex $u=[1,1,0,0]$. Since the diameter of $A(q)$ is $3$, it is sufficient to show that for every path $P$ of length at most $3$ in $A(q)$ from $u$ to an arbitrary vertex $w\in O_5$ and passing through a vertex outside $O_5$ the distance between $u$ and $w$ is still at most $3$ in $C(G,S^*)$. By our earlier examination of possibilities, and referring also to the notation introduced in the previous paragraph, such a path $P$ can only have one of the following eight forms:

1) $P=uu_1w$, with $w\in L(u_1)$,

2) $P=uu_2w$, with $w\in L(u_2)$,

3) $P=uu_1w'w$ for some $w'\in L(u_1)$, with $w\notin L(u_1)\cup\{u\}$

4) $P=uu_2w'w$ for some $w'\in L(u_2)$, with $w\notin L(u_2)\cup\{u\}$

5) $P=uu_1u'_2w$, with $w\in L(u'_2)$,

6) $P=uu_2u'_1w$, with $w\in L(u'_1)$,

7) $P=uvv_1w$ for some $v\in O_5$, with $w\in L(v_1)$, or

8) $P=uvv_2w$ for some $v\in O_5$, with $w\in L(v_2)$.
\smallskip

Since $S_1,S_2\subset S^*$, we may apply Lemma \ref{stab} to conclude that the subgraphs of $C(G,S^*)$ induced by the subsets $L(u_1)\cup\{u\}$ and $L(u_2)\cup\{u\}$ are isomorphic to the Cayley graphs $C(H_1,S_1)$ and $C(H_2,S_2)$, both of diameter $2$. This implies that the distance between $u$ and $w$ is at most $3$ in the cases 1) -- 4). Similarly, since $H_2$ is also equal to $Stab_G(u'_1)$ and $H_1$ is equal to $Stab_G(u'_2)$, by Lemma \ref{stab} (which also describes the action of $H_2$ and $H_1$) the distances in $C(G,S^*)$ between vertices in the sets $L(u'_1)$ and $L(u'_2)$ are at most $2$. Having included the elements $M_1=M(1;1,1)$ and $M_2=M(1;1,0)$ in $S^*$, we have an edge from $u=[1,1,0,0]$ to $uM_1=[1,0,1,1]\in L(u'_1)$ and an edge from $u$ to $uM_2=[1,0,1,0]\in L(u'_2)$. It follows that we have paths of length $3$ from $u$ to $w$ inside $C(G,S^*)$ also in the cases 5) and 6).
\smallskip

Let $v$ be a vertex as in the cases 7) and 8), incident with $u$ in $B(q)$. By the regular action of $G$ on the vertex set of $B(q)$ there exists a (unique) $g\in G$ such that $v=ug$. From the previous paragraph we know that, for $i=1,2$, the subgraph of $C(G,S^*)$ induced by the vertex set $L(u_i)\cup\{u\}$ is isomorphic to the Cayley graph $C(H_i,S_i)$. But since $L(v_i)\cup\{v\}=(L(u_i)\cup\{u\})g$ for $i=1,2$, the subgraph of $C(G,S^*)$ induced by the set $L(v_i)\cup\{v\}$ is isomorphic to the one induced by the set $L(u_i)\cup\{u\}$, which is the Cayley graph $C(H_i,S_i)$. It follows that in the cases 7) and 8), for $i=1,2$, the $vv_iw$-part of the path $P$ can be replaced by a path of length at most two in the subgraph of $C(G,S^*)$ induced by the set $L(v_i)\cup\{v\}$, since the diameter of $C(H_i,S_i)$ is $2$. This completes the proof. \ebox

Letting $d=q+2c_{q-1}+3\le q + 4\lceil\sqrt{q-1}\rceil +3$ for $q=2^{2n+1}$ and asymptotically expressing the order $q^2(q-1)$ of $C(G,S^*)$ in terms of $d$, we obtain:

\bc\label{cor1} There exists an infinite increasing sequence of values of $d$ for which there are Cayley graphs of degree $d$, diameter $3$ and order $d^3-O(d^{2.5})$ as $d\to\infty$. \hfill $\Box$
\ec

We thus have ${\rm lim\ sup}_{d\to\infty}\ Cay(d,3)/M(d,3)=1$, which means that the Moore bound for diameter $3$ can be asymptotically approached by Cayley graphs.
\smallskip

\section{Remarks}

Our construction of an infinite sequence of Cayley graphs of diameter 3, degree $q+o(q)$ and order $q^3-o(q^3)$ was based on the existence of a group $G$ of order $q^2(q-1)$ for $q=2^{2n+1}$, regular on the vertex set of a subgraph obtained from the point-line incidence graph of a generalised quadrangle $W(q)$ by factorisation by a polarity. Moreover, the group $G$ is (isomorphic to) a subgroup of index $q^2+1$ of the Suzuki group ${\rm Sz}(q)={^2}B_2(q)$, a simple group of order $q^2(q^2+1)(q-1)$, acting doubly transitively on the Suzuki-Tits ovoid (of order $q^2+1$) in the generalised quadrangle $W(q)$; see \cite{M} for more details.
\smallskip

It is tempting to consider an analogous approach for constructing an infinite sequence of Cayley graphs of diameter $5$, degree $q+o(q)$ and order $q^5-o(q^5)$ from a suitable group acting regularly on some subgraph obtained from the incidence graph $I(q)$ of a {\sl generalised hexagon} $H(q)$ factored by a polarity. By known results summarised in \cite{M}, a generalised hexagon $H(q)$ admits a polarity if and only if $q=3^{2n+1}$, $n\ge 1$. For such $q$, the graph $I(q)$, of diameter $5$ and order $q^5+\ldots +q+1$, has $q^5+q^4+q^2+q$ vertices of degree $q+1$ together with the $q^3+1$ vertices of degree $q$ that  correspond to the Ree-Tits ovoid formed by the absolute (that is, self-polar) points of $H(q)$. It turns out that, up to field automorphisms, the automorphism group of $I(q)$ can be identified with the collineation group that preserves the Ree-Tits ovoid, which is the Ree group ${\rm Re}(q)={^2}G_2(q)$.
\smallskip

Unfortunately, by the classification of maximal subgroups of the Ree groups \cite{LN,K}, the group ${\rm Re}(q)={^2}G_2(q)$, a simple group of order $q^3(q^3+1)(q-1)$, does not contain a subgroup of order $O(q^5)$ for $q\to\infty$. The approach that works for construction of large Cayley graphs of diameter $2$ from generalised triangles in \cite{BS}, and of diameter $3$ from generalised quadrangles in this article, thus does not carry over to an analogous construction for diameter $5$ from generalised hexagons.

\bigskip
\bigskip

\noindent{\bf Acknowledgement.}~~ The second and the third authors acknowledge support from the APVV Research Grant 0136-12 and the VEGA Research Grants 1/0065/13 and 1/0007/14.

\end{document}